\documentclass[12pt,thmsa]{article}%
\usepackage{amsfonts}
\usepackage{sw20bams}
\usepackage{amsmath}
\usepackage{amssymb}
\usepackage{graphicx}%
\providecommand{\U}[1]{\protect\rule{.1in}{.1in}}
\begin{document}

\author{Steven R. Finch}
\title{The Logarithmic Spiral Conjecture}
\date{March 12, 2016}
\maketitle

\begin{abstract}
When searching for a planar line, if given no further information, one should
adopt a logarithmic spiral strategy (although unproven).

\end{abstract}

\footnotetext{Copyright \copyright \ 2005, 2016 by Steven R. Finch. All rights
reserved.}This brief paper is concerned entirely with geometry in the plane
and continues a thought in \cite{FW}. If a line intersects a circle in one or
two points, we say that the line \textbf{strikes} the circle. If a line
intersects a circle in exactly one point (that is, if the line is tangent to
the circle), we say that the line \textbf{touches} the circle.

Let $f$ be a nonnegative, continuously differentiable function on $\mathbb{R}$
satisfying
\[%
\begin{array}
[c]{ccc}%
\lim\limits_{\theta\rightarrow-\infty}f(\theta)=0, &  & \lim\limits_{\theta
\rightarrow\infty}f(\theta)=\infty.
\end{array}
\]
The polar curve $r=f(\theta)$ intersects every line in the plane, that is, $f$
is a \textbf{spiral}. (Reason: for each $R>0$, there exists $\Theta$ so large
that $\theta>\Theta$ implies $f(\theta)>R$. Any line striking the circle $r=R$
must therefore intersect the curve $r=f(\theta)$. Since $R$ was arbitrary, the
statement follows.) Existence of intersection points is only the beginning of
our study.

Consider the set $\Sigma$ of all lines that strike the circle $r=R$. The
spiral $r=f(\theta)$ possesses a first intersection point $\theta$ with each
line in $\Sigma$; let $\theta_{1}$ denote the supremum of all such $\theta$.
Loosely put, $\theta_{1}$ constitutes the worst case scenario when seeking all
members of $\Sigma$ via the search strategy $r=f(\theta)$. Clearly $\theta
_{1}$ depends on $R$ and $\theta_{1}=-\infty$ when $R=0$.

The cost of finding all lines in $\Sigma$, starting from the origin, can be
quantified by the arclength
\[
\Lambda(f)=%
{\displaystyle\int\limits_{-\infty}^{\theta_{1}}}
\sqrt{f(\theta)^{2}+f^{\prime}(\theta)^{2}}\,d\theta.
\]
We naturally wish to minimize $\Lambda(f)$ as a function of $f$, for fixed
$R$. Our focus is on the following asymptotic inequality.\medskip

\noindent\textbf{Conjecture 1.} \
\[
\lim_{R\rightarrow\infty}\frac{\Lambda(f)}{R}\geq13.8111351795...
\]
\textit{with equality if and only if} $f(\theta)\sim Ce^{\kappa\theta}$
\textit{as} $\theta\rightarrow\infty$\textit{, where} $\kappa=0.2124695594...$
\textit{and} $C>0$\textit{ is arbitrary}$.$\medskip

The two numerical constants appear precisely in \cite{Fin}, along with
detailed treatment of the special case of a logarithmic spiral $f(\theta
)=e^{\kappa\theta}$. Difficulties arise in the general case, owing to the vast
variety of spirals permitted

A sketch of a geometric proof of Conjecture 1 was published in \cite{BY1, BY2,
BY3}. The first part claimed that an optimal spiral must be similar with
respect to both rotations and dilations about the origin; the second part
claimed that such a highly symmetric spiral must necessarily be a logarithmic
spiral. The second part, in fact, is true via the solution of a well-known
functional equation \cite{You}. We doubt, however, that any purely geometric
proof of the first part can be rigorously correct (although appealing). A more
careful analysis, based on the calculus of variations, is perhaps mandatory.

\subsection{Examples}

We repeat certain steps employed in \cite{Fin}, suitably generalized.\medskip

\noindent\textbf{Lemma 2.} \textit{The distance between the line} $Ax+By+C=0$
\textit{and the origin is} $|C|/\sqrt{A^{2}+B^{2}}$.\medskip

\noindent\textbf{Lemma 3.} \textit{The equation of a line tangent to the
spiral} $r=f(\theta)$ \textit{is} $y-f(\theta)\sin(\theta)=m(x-f(\theta
)\cos(\theta))$\textit{, where }$\theta$\textit{ corresponds to the point of
tangency and the slope is given by}
\[
m=\frac{f^{\prime}(\theta)\sin(\theta)+f(\theta)\cos(\theta)}{f^{\prime
}(\theta)\cos(\theta)-f(\theta)\sin(\theta)}.
\]

\noindent\textbf{Proof of Lemma 3.} Clearly
\[
\frac{dy}{dx}=\frac{dy/d\theta}{dx/d\theta}=\frac{(f(\theta)\sin
(\theta))^{\prime}}{(f(\theta)\cos(\theta))^{\prime}}=\frac{f^{\prime}%
(\theta)\sin(\theta)+f(\theta)\cos(\theta)}{f^{\prime}(\theta)\cos
(\theta)-f(\theta)\sin(\theta)}.\medskip
\]
\noindent\textbf{Theorem 4}.\textit{ Let} $L$ \textit{denote the first line
that is both tangent to the spiral} $r=f(\theta)$ \textit{and tangent to the
circle} $r=R$. \textit{The tangency point} $\theta_{0}$ \textit{of} $L$
\textit{with the spiral satisfies the equation}
\[
R^{2}(f(\theta)^{2}+f^{\prime}(\theta)^{2})=f(\theta)^{4}.
\]

\noindent\textbf{Proof of Theorem 4.} Apply Lemma 2 with $A=m$, $B=-1$ and
$C=f(\theta)(\sin(\theta)-m\cos(\theta))$ to obtain $(1+m^{2})R^{2}%
=f(\theta)^{2}(\sin(\theta)-m\cos(\theta))^{2}$. Substituting the expression
for $m$ from Lemma 3 gives the desired equation.\newline

We emphasize that, on the one hand, $\theta_{0}$ is where the spiral first
intersects a line that touches the circle $r=R$ (the touching occurs
elsewhere). On the other hand, $\theta_{1}$ is just above where the spiral
last intersects a new line that strikes the circle $r=R$ (the striking, again,
occurs elsewhere). If the function $f$ is strictly increasing, then in the
interval $\theta_{0}<\theta<\theta_{1}$, the spiral intersects all other lines
that touch $r=R$; at $\theta=\theta_{1}$, repetition begins so we stop there.
Suppose that we are given a spiral $r=f(\theta)$ for which $f(\theta
)\not \sim Ce^{\kappa\theta}$ for any $\kappa>0$, $C>0$. Clearly
\[
\frac{\Lambda(f)}R\geq\frac1R%
{\displaystyle\int\limits_{-\infty}^{\theta_{0}}}
\sqrt{f(\theta)^{2}+f^{\prime}(\theta)^{2}}\,d\theta,
\]
and thus if we demonstrate that the right hand side $\rightarrow\infty$ or is
at least $>13.82$, then this is consistent with Conjecture 1.

As a first example, consider Archimedes' spiral
\[
f(\theta)=\left\{
\begin{array}
[c]{lll}%
\kappa\theta &  & \text{if }\theta\geq0,\\
0 &  & \text{if }\theta<0.
\end{array}
\right.
\]
From Theorem 4, it follows that $R^{2}(1+\theta^{2})=\kappa^{2}\theta^{4}$ and
hence
\[
\theta_{0}=\frac{R}{\kappa}\sqrt{\frac{1}{2}\left(  1+\sqrt{1+\frac
{4\kappa^{2}}{R^{2}}}\right)  }\geq\frac{R}{\kappa}\sqrt{\frac{1}{2}\left(
1+1\right)  }\geq\frac{R}{\kappa}.
\]
Consequently, the normalized arclength is bounded from below by
\[
\frac{\kappa}{R}%
{\displaystyle\int\limits_{0}^{\theta_{0}}}
\sqrt{1+\theta^{2}}\,d\theta\geq\frac{\kappa}{R}%
{\displaystyle\int\limits_{0}^{\theta_{0}}}
\theta\,d\theta=\frac{\kappa}{2R}\theta_{0}^{2}\geq\frac{R}{2\kappa
}\rightarrow\infty
\]
as $R\rightarrow\infty$. Alternatively, we can avoid solving for $\theta_{0}$
altogether:\ From $R^{2}(1+\theta^{2})=\kappa^{2}\theta^{4}$, deduce that
\[
R=\frac{\kappa\,\theta^{2}}{\sqrt{1+\theta^{2}}}\leq\kappa\,\theta^{2}%
\]
and hence that $\theta\rightarrow\infty$ as $R\rightarrow\infty$. Here we
obtain
\[
\frac{\kappa}{R}%
{\displaystyle\int\limits_{0}^{\theta_{0}}}
\theta\,d\theta=\frac{\kappa}{2R}\theta_{0}^{2}=\frac{\kappa}{2}\frac
{\sqrt{1+\theta_{0}^{2}}}{\kappa\,\theta_{0}^{2}}\theta_{0}^{2}=\frac{1}%
{2}\sqrt{1+\theta_{0}^{2}}\geq\frac{\theta_{0}}{2}\rightarrow\infty
\]
as $\theta_{0}\rightarrow\infty$ (and thus as $R\rightarrow\infty$). This
latter device will be useful in the following examples. See Figure 1 for an
illustration.%
\begin{figure}[ptb]%
\centering
\includegraphics[
height=5.3696in,
width=5.2978in
]%
{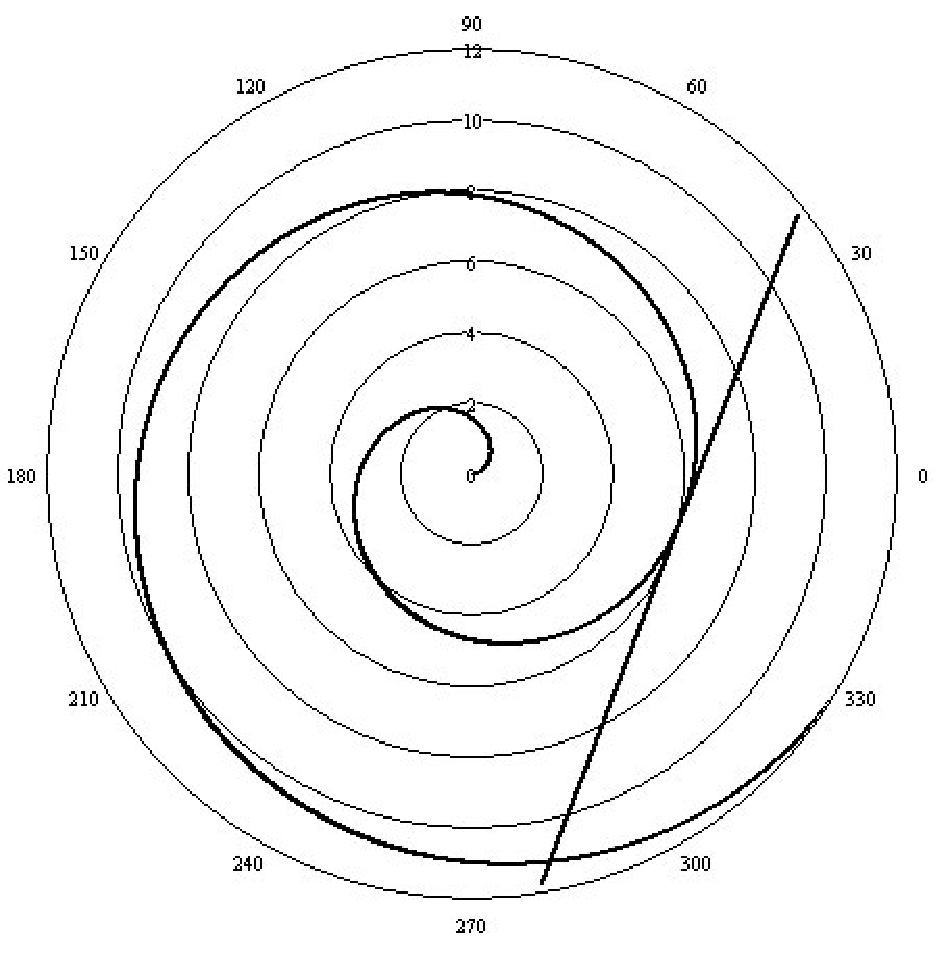}%
\caption{The first contact point that the spiral $r=\theta$ has with a line
tangent to the circle $R=6$ is at $\theta_{0}=348.4^{\circ}$. The second
contact point with the line is at $\theta_{1}=641.5^{\circ}$. Incidently, the
line is tangent to $R=6$ at $339.1^{\circ}<\theta_{0}$.}%
\end{figure}

Consider next the spiral
\[
f(\theta)=\left\{
\begin{array}
[c]{lll}%
e^{\theta^{a}} &  & \text{if }\theta\geq0,\\
e^{-\left|  \theta\right|  ^{a}} &  & \text{if }\theta<0
\end{array}
\right.
\]
for a fixed exponent $a>0$. From Theorem 4, it follows that $R^{2}%
(1+a^{2}\theta^{2a-2})=e^{2\theta^{a}}$, that is,
\[
R=\frac{e^{\theta^{a}}}{\sqrt{1+a^{2}\theta^{2a-2}}}\leq e^{\theta^{a}}.
\]
Hence $\theta\rightarrow\infty$ as $R\rightarrow\infty$. If $a>1$, the
normalized arclength is bounded from below by
\begin{align*}
&  \frac1R\left(
{\displaystyle\int\limits_{-\infty}^{0}}
\sqrt{1+a^{2}\left|  \theta\right|  ^{2a-2}}\,e^{-\left|  \theta\right|  ^{a}%
}d\theta+%
{\displaystyle\int\limits_{0}^{\theta_{0}}}
\sqrt{1+a^{2}\theta^{2a-2}}\,e^{\theta^{a}}d\theta\right) \\
&  \geq\frac1R\left(
{\displaystyle\int\limits_{-\infty}^{0}}
a\,\left|  \theta\right|  ^{a-1}e^{-\left|  \theta\right|  ^{a}}d\theta+%
{\displaystyle\int\limits_{0}^{\theta_{0}}}
a\,\theta^{a-1}e^{\theta^{a}}d\theta\right) \\
&  =\frac1R\left(  1+e^{\theta_{0}^{a}}-1\right)  =\sqrt{1+a^{2}\theta
_{0}^{2a-2}}\geq a\,\theta_{0}^{a-1}\rightarrow\infty
\end{align*}
as $\theta_{0}\rightarrow\infty$ (and thus as $R\rightarrow\infty$). If
$0<a<1$, the normalized arclength is bounded by
\begin{align*}
&  \ \frac1R\left(
{\displaystyle\int\limits_{-\infty}^{0}}
\sqrt{1+a^{2}\left|  \theta\right|  ^{2a-2}}\,e^{-\left|  \theta\right|  ^{a}%
}d\theta+%
{\displaystyle\int\limits_{0}^{\theta_{0}}}
\sqrt{1+a^{2}\theta^{2a-2}}\,e^{\theta^{a}}d\theta\right) \\
\  &  \geq\frac1R\left(
{\displaystyle\int\limits_{-\infty}^{0}}
e^{-\left|  \theta\right|  ^{a}}d\theta+%
{\displaystyle\int\limits_{0}^{\theta_{0}}}
e^{\theta^{a}}d\theta\right)  \geq\frac1R\left(  0+%
{\displaystyle\int\limits_{0}^{\theta_{0}}}
e^{\theta^{a}}d\theta\right)
\end{align*}
and we have asymptotics
\[
\frac1R%
{\displaystyle\int\limits_{0}^{\theta_{0}}}
e^{\theta^{a}}d\theta\sim\frac1R\left(  \frac1a\theta_{0}^{1-a}e^{\theta
_{0}^{a}}\right)  =\frac1a\sqrt{1+a^{2}\theta_{0}^{2a-2}}\,\theta_{0}%
^{1-a}\rightarrow\infty
\]
as $\theta_{0}\rightarrow\infty$ (and thus as $R\rightarrow\infty$). Only the
case $a=1$ remains, which is covered in \cite{Fin}. This is compelling (but
not completely convincing) evidence that the Logarithmic Spiral Conjecture is valid.

Consider finally the spiral
\[
f(\theta)=\left\{
\begin{array}
[c]{lll}%
\theta^{b}e^{\theta} &  & \text{if }\theta\geq0,\\
0 &  & \text{if }\theta<0
\end{array}
\right.
\]
for a fixed exponent $b>0$. From Theorem 4, we obtain
\[
R=\frac{\theta^{b+1}e^{\theta}}{\sqrt{\theta^{2}+(b+\theta)^{2}}}\leq
\frac1{\sqrt{2}}\theta^{b}e^{\theta},
\]
hence $\theta\rightarrow\infty$ as $R\rightarrow\infty$. Clearly $\theta
_{1}\geq\theta_{0}+\pi$ on geometric grounds. Therefore the normalized
arclength is bounded from below by
\[
\frac1R%
{\displaystyle\int\limits_{0}^{\theta_{0}+\pi}}
\sqrt{\theta^{2}+(b+\theta)^{2}}\,\theta^{b-1}e^{\theta}d\theta\geq\frac
{\sqrt{2}}R%
{\displaystyle\int\limits_{0}^{\theta_{0}+\pi}}
\theta^{b}e^{\theta}d\theta
\]
and we have asymptotics
\begin{align*}
\frac{\sqrt{2}}R%
{\displaystyle\int\limits_{0}^{\theta_{0}+\pi}}
\theta^{b}e^{\theta}d\theta &  \sim\frac{\sqrt{2}}R\left(  (\theta_{0}%
+\pi)^{b}e^{\theta_{0}+\pi}\right) \\
\  &  =\sqrt{2}\,\frac{(\theta_{0}+\pi)^{b}e^{\theta_{0}+\pi}}{\theta
_{0}^{b+1}e^{\theta_{0}}}\sqrt{\theta_{0}^{2}+(b+\theta_{0})^{2}}\\
\  &  \rightarrow2e^{\pi}>13.82
\end{align*}
as $\theta_{0}\rightarrow\infty$.

The logarithmic spiral appears with regard to another planar search problem
\cite{GC}, but the techniques of Gal \&\ Chazan do not seem to apply here. A
min-mean analog of our Conjecture 1 could also be formulated, starting with
\cite{Fin}.

\subsection{Acknowledgements}

I am grateful to Ricardo Baeza-Yates, Li-Yan Zhu and John Shonder for their
assistance. A discussion during my Visiting Lecture at Oberlin College in
November 2004 was also very helpful.

\end{document}